\documentclass[11pt,twoside]{article}

\usepackage{amsmath}
\usepackage{amsthm}
\usepackage{amsfonts, amssymb}
\usepackage{mathrsfs}
\usepackage[all]{xy}

\def\titulo#1{\noindent{\bf\LARGE{#1}} \bigskip \thispagestyle{plain}}
\def\autor#1{\noindent{\sc #1}\smallskip}
\def\direccion#1{\noindent #1\bigskip}
\def\email#1{\vspace{1cm}\noindent E-mail address: {\sf #1} \bigskip}

\usepackage{latexsym}
\usepackage[dvips]{graphicx}

\theoremstyle{plain}
\newtheorem{lema}{Lemma}[section]
\newtheorem{prop}[lema]{Proposition}
\newtheorem{teo}[lema]{Theorem}
\newtheorem{coro}[lema]{Corollary}
\theoremstyle{remark}

\newtheorem{obs}[lema]{Remark}
\newtheorem{obsi}[lema]{Important Remark}

\theoremstyle{definition}
\newtheorem{defi}[lema]{Definition}
\newtheorem{ej}[lema]{Example}

\def\o{\mathfrak{o}}
\def\e{\mathfrak{e}}

\def\k{\mathcal{K}}
\def\x{\mathcal{X}}
\def\etft0{\mathnormal{etfT_0}}

\def\S{\mathbb{S}}

\begin{document}

\titulo{Minimal Finite Models}

\autor{Jonathan Ariel Barmak, Elias Gabriel Minian}

\direccion{Departamento  de Matem\'atica.\\
 FCEyN, Universidad de Buenos Aires. \\ Buenos
Aires, Argentina}

\begin{abstract}
\noindent We characterize the smallest finite spaces with the same homotopy groups of the spheres. Similarly, we describe the minimal finite models of any 
finite graph.
We also develop new combinatorial techniques based on finite spaces to study classical invariants of general topological spaces. 
\end{abstract}

\noindent{\small \it 2000 Mathematics Subject Classification.
 \rm 55P10, 55P15, 18B35, 18G30.}

\noindent{\small \it Key words and phrases. \rm Finite Spaces, Weak Homotopy Types, Spheres, Graphs, Posets.}

\section{Introduction}

This paper deals with finite topological spaces and their application to the homotopy theory of (general) topological spaces. One of its main goals is to
characterize the  \it minimal finite models \rm of some spaces such as spheres and finite graphs. A minimal finite model of a topological space $Y$ is a finite space with 
the smallest number of points that 
is weakly homotopy equivalent to $Y$.

It is well known \cite{May2, Mcc} that finite spaces are  related to simplicial complexes. Explicitly, following McCord \cite{Mcc} 
one can associate to any finite space $X$ a finite simplicial complex $\k(X)$ and a weak homotopy equivalence $|\k(X)|\to X$. Moreover,
given a finite simplicial complex $K$, there is a finite space $\x(K)$ and a weak homotopy equivalence $|K|\to \x(K)$.

One can find this way \it finite models \rm of topological spaces, i.e. finite topological spaces with the same weak homotopy type.
McCord also exhibits in \cite{Mcc} finite models for the spheres $S^n$, denoted $\S^n S^0$, with only $2n+2$ points. 

\medskip

In his  series of notes on 
finite spaces \cite{May,May2,May3}, J.P. May conjectures that $\S^n S^0$ is, in our terminology, a minimal finite model for  the $n$-dimensional
sphere. We prove that this conjecture is true. In fact, we prove the following stronger result:

\begin{teo} \label{modeloesfera} Any space with  the same homotopy groups of $S^n$ has at least $2n+2$ points. Moreover, $\S^n S^0$ is the unique space with $2n+2$ points with
this
property.
\end{teo}

In particular, $\S^n S^0$ is a minimal finite model of $S^n$ and it is unique.

\medskip

We also obtain a similar result for minimal finite models of finite graphs. It is well known that finite graphs are homotopy equivalent to a wedge sum of finitely
many copies of one-dimensional spheres $\bigvee\limits_{i=1}^{m}S^1$ (cf. \cite{May4}). Their finite minimal models are characterized as follows.

\begin{teo}\label{modelografo}
Let $n\in \mathbb{N}$. A finite $T_0$-space $X$ is a minimal finite model of $\bigvee\limits_{i=1}^{n}S^1$ if and only if $h(X)=2$, $\# X=min\{ i+j \ | 
\ (i-1)(j-1)\ge n\}$ and $\# \mathtt{E} (\mathcal{H}(X))=\# X +n-1$.
\end{teo}

Here $h(X)$ denotes the height of $X$ (viewed as a poset), $\# X$ denotes the number of points of $X$ and $\# \mathtt{E}(\mathcal{H}(X))$ the number of edges of the
Hasse diagram of $X$. In particular, one can explicitly compute the number of points of any minimal finite model of a finite graph.

\medskip

Note that, in general, minimal finite models are not unique. If $X$ is a finite model of a space, then so is $X^{op}$, which is the opposite preorder of $X$ (see \ref{op1} and \ref{op2}). Moreover, a space can have more than two minimal finite models. To illustrate this, we exhibit in the last section of this paper the three minimal models 
of $\bigvee\limits_{i=1}^{3}S^1$, each of which with 6 points and 8 edges.

\bigskip

The main reason of investigating finite models of spaces with the same weak homotopy type instead of finite models with the same homotopy type is that
the homotopy type of finite spaces rarely occurs in general spaces. More precisely, we prove in section 2 the following result:

\begin{teo} If $X$ is a $T_1$, connected and non contractible space, then it does not have the homotopy type of any finite space.
\end{teo}

In particular, finite spaces do not have the same homotopy type of any connected non contractible CW-complex.

\medskip

We have to make an important  remark about a paper   of Osaki \cite{Osa}. In that article he introduces two methods of reduction which allow one to shrink 
a finite $T_0$-space to a smaller weakly equivalent space and he asks whether, by a sequence of these kinds of reductions, each finite $T_0$-space 
can be reduced to the smallest one with the same homotopy groups. In section 2 of this paper, we exhibit an example which shows that the answer to
his question is negative. Therefore, his methods of reduction are not always effective and we could not apply them to prove Theorems
\ref{modeloesfera} and \ref{modelografo} of above.

\medskip

Finally, we would like to emphasize that the methods and tools that we develop here are, in our opinion, as important as the results that we mentioned 
above. Algebraic topologists are  commonly used to work with the characteristic combinatorial methods  of simplicial complexes and cellular spaces. 
We think that the combinatorial tools that one can deduce from finite spaces can be, in many situations,  even more efficient for investigating 
 homotopy and homology theory of general topological spaces.

\medskip

To illustrate this, consider the following result, proved in section 4, which is one of the key points in the solution of the problem of the 
minimal finite 
models of graphs: 

\begin{prop}
Let $X$ be a connected finite $T_0$-space and let $x_0\neq x\in X$ such that $x$ is neither maximal nor minimal in $X$. 
Then the inclusion map of the associated simplicial complexes $\k (X\smallsetminus \{x\})\subseteq \k (X)$ induces an epimorphism 
 $i_*:E(\k(X\smallsetminus \{x\}),x_0)\to E(\k(X),x_0)$ between their edge-path (fundamental) groups. 
\end{prop}

The conditions of maximality or minimality of points in a finite space, as well as the notion of \it beat point \rm introduced by Stong \cite{Sto}, are
hard to express in terms of simplicial complexes. In this direction, we will show in section 3 how to compute combinatorially 
the fundamental group of a finite $T_0$-space from its Hasse diagram.

\section{Preliminaries and the problem of the spheres}

Let $X$ be a finite space. For each $x\in  X$ we denote by  $U_x$ the minimal open set of $x$, defined as the 
intersection of all open sets containing $x$ (cf. \cite{May,Mcc}).

Given a topology in a finite set $X$, the associated preorder in $X$ is defined by $x\le y$ if $x\in U_y$. 
Alexandroff \cite{Ale} proved that this aplication is a one to one correspondence between topologies and preorders in $X$.
Moreover, $T_0$-topologies correspond to (partial) orders. 

Therefore we will regard finite spaces as finite preorders and viceversa.

Note that a function $f:X\to Y$ between finite spaces is continuous if and only if it is order preserving.

\begin{defi} \label{op1}
Let $X$ be a finite space. We define $X^{op}$ as the space whose underlying set is  $X$ but with the  opposite preorder.
\end{defi}

We recall that the Hasse diagram of a finite poset $P$ is the digraph whose set of vertices is $P$ and whose 
edges are the ordered pairs $(x,y)\in P$ with $x<y$ such that there is no $z\in P$ with $x<z<y$. 
We define $\mathcal{H}(X)=(\mathtt{V}(\mathcal{H}(X)),\mathtt{E}(\mathcal{H}(X)))$, the \textit{Hasse diagram} of a finite $T_0$-space $X$, as the Hasse diagram of the associated order of $X$.

In order to indicate the orientation of an edge of $\mathcal{H}(X)$ in a figure, we will put $y$ over $x$ if $(x,y)\in \mathtt{E}(\mathcal{H}(X))$.

\begin{ej}
Let $X=\{a,b,c,d\}$ be the space whose open sets are $\emptyset$, $\{a,b,c,d\}$, $\{b,d\}$, $\{c\}$, $\{d\}$, $\{b,c,d\}$ and $\{c,d\}$. The Hasse diagram of $X$ is 

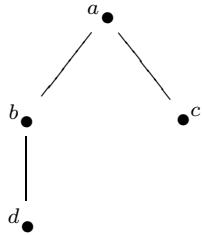
\begin{displaymath}
\xymatrix@C=10pt{ & ^a \bullet \ \ar@{-}[rd] \ar@{-}[ld] \\ ^b \bullet \ \ar@{-}[d] & & \ \bullet ^c \\ ^d \bullet \ }
\end{displaymath}

\end{ej}

Following McCord \cite{Mcc}, one can see that, in order to investigate homotopy types of finite spaces, it suffices to study $T_0$-spaces. 
Stong developed in \cite{Sto} a very powerful tool to classify the homotopy types of $T_0$-spaces. We recall from \cite{Sto} and May's notes
\cite{May}  the following definitions and results:

\begin{defi} (Stong)
Let $X$ be a finite $T_0$-space. 

A point $x\in X$ is called an \textit{up beat point} if there exists $y\in X$, $y>x$ such that $z> x$ implies $z\ge y$.

Analogously, a point $x\in X$ will be called a \textit{down beat point} if there exists $y\in X$, $y<x$ such that $z< x$ implies $z\le y$.

A finite $T_0$-space $X$ is called a \textit{minimal finite space} if it has no beat points.
\end{defi}

For any finite space $X$, Stong defines its core as a minimal finite space 
which is a strong deformation retract of $X$. He shows that any finite space has a core, and that the only  map $f:X\rightarrow X$ 
between minimal finite spaces which is homotopic to the identity is the identity itself. Then, the core $X_c$ of a space $X$, is unique up to homeomorphism and it is the space of minimum cardinality that is homotopy equivalent to $X$.

For example, we obtain that a finite space is contractible if and only if its core is a point. Since the core of a finite space is the disjoint union of the cores of its connected components, we can deduce the following 

\begin{lema}
Let $X$ be a finite space such that $X_c$ is discrete. Then $X$ is a disjoint union of contractible spaces. 
\end{lema}

As we pointed out in the introduction, finite spaces do not have in general the same homotopy type of $T_1$-spaces:

\begin{teo}
Let $X$ be a finite space and let $Y$ be a $T_1$-space homotopy equivalent to $X$. Then $X$ is a disjoin union of contractible spaces. 
\end{teo}

\begin{proof}
Since $X\simeq Y$, then $X_c\simeq Y$. Let $f:X_c\rightarrow Y$ be a homotopy equivalence with homotopy inverse $g$. 
Then $gf=1_{X_c}$ since $X_c$ is a minimal finite space. 
Since $f$ is a one to one map from  $X_c$ to a $T_1$-space, it follows that $X_c$ is also $T_1$ and therefore discrete. 
Now the result follows from the previous lemma.
\end{proof}

\begin{coro}
Let $Y$ be a  connected and non contractible $T_1$-space. Then $Y$ does not have the same homotopy type of any finite space.
\end{coro}

\begin{proof}
Follows immediately from the previous Theorem.
\end{proof}

For example, for any $n\ge 1$, the $n$-dimensional sphere $S^n$ does not have the homotopy type of any finite space. Although, $S^n$ 
does have, as any finite polyhedron,  the same weak homotopy type of some finite space.

\begin{defi}
Let $X$ be a space. We say that a finite space $Y$ is a \textit{finite model} of $X$ if it is weakly equivalent to $X$.

We say that $Y$ is a \textit{minimal finite model} if it is a finite model of minimum cardinality.
\end{defi}

By weakly (homotopy) equivalent, we mean a topological space $Y$ such that there is a finite sequence $X=X_0,X_1,\ldots,X_r=Y$ and weak homotopy
equivalences $X_i\to X_{i+1}$ or $X_{i+1}\to X_i$ for each $i=0,\ldots, r-1$.

\medskip

For example, the singleton is the unique minimal finite model of every contractible space.

Since every finite space is homotopy equivalent to its core, which is a smaller space, we have the following

\begin{prop} \label{mfmimplicaefm}
Every minimal finite model is a minimal finite space.
\end{prop}

In  \cite{Mcc} McCord associates to each finite $T_0$-space $X$ a simplicial complex $\k (X)$ whose simplices are the non empty chains of $X$ 
and proves that $|\k (X)|$ is weakly equivalent to $X$.\\

Since $\k (X)=\k (X^{op})$ we get the following

\begin{prop} \label{op2}
If $X$ is a minimal finite model of a space $Y$, then so is $X^{op}$.
\end{prop}

\begin{ej}
The 5-point $T_0$-space $X$, whose Hasse diagram is 
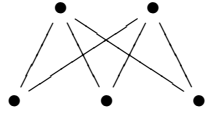
\begin{displaymath}
\xymatrix@C=6pt{ & \bullet \ar@{-}[dl] \ar@{-}[dr] \ar@{-}[drrr] & & \bullet \ar@{-}[dlll] \ar@{-}[dl] \ar@{-}[dr]\\ 
		\bullet & & \bullet & & \bullet} 
\end{displaymath}
has an associated polyhedron  $|\k (X)|$, which is homotopy equivalent to $S^1\vee S^1$. Therefore, $X$ is a finite model of $S^1\vee S^1$. 
In fact,  it is a minimal finite model since every space with less than 5 points is either contractible, or non connected or weakly equivalent 
to $S^1$. However, this minimal finite model is  not unique since $X^{op}$ is another minimal finite model  not homeomorphic to $X$. \\

We will generalize this result later, when we characterize the minimal finite models of graphs.
\end{ej}

Note that, by Whitehead Theorem, if  $X$ is a  finite model of a (compact) polyhedron $Y$, then $Y$ is homotopy equivalent to $|\k (X)|$.

\bigskip

Let $X$ be a finite space. The \textit{non-Hausdorff suspension} $\mathbb{S}X$ of $X$ is  the finite space $X\cup \{+,-\}$ 
whose open sets are those of $X$ together with $X\cup \{+\}$, $X\cup \{-\}$ and $X\cup \{+,-\}$.

The non-Hausdorff suspension of order $n$ is defined recursively by $\mathbb{S}^{n}X=\mathbb{S}(\mathbb{S}^{n-1}X)$.


\medskip

In \cite{Mcc}, McCord proved that the $(2n+2)$-point space $\mathbb{S}^nS^0$ is a finite model of $S^n$.
In \cite{May2} May conjectures that $\mathbb{S}^nS^0$ is a minimal finite model of the sphere. We will show that this conjecture is true. In fact, we 
prove a stronger result. Namely, we will see that
any space with the same homotopy groups of $S^n$ has at least $2n+2$ points. Moreover, if it has exactly $2n+2$ points then it has to be homeomorphic to
$\S^n S^0$.

Before we proceed with the proof of the conjecture, we would like to make some remarks about Osaki's methods of reduction \cite{Osa}.

In \cite{Osa} Osaki proves  the following results.

\begin{teo} (Osaki)
Let $X$ be a finite $T_0$-space. Suppose there exists  $x\in X$ such that $\pi_n(U_x\cap U_y)=0$ for all $n\ge 0$ and for all $y\in X$. 
Then the quotient map $p:X\rightarrow X/U_x$ is a weak homotopy equivalence. 
\end{teo}

\begin{teo} (Osaki)
Let $X$ be a finite $T_0$-space. Suppose there exists $x\in X$ such that $\pi_n(\overline{\{x\}}\cap \overline{\{y\}})=0$ for all $n\ge 0$ and for all $y\in X$. 
Then the quotient map $p:X\rightarrow X/\overline{\{x\}}$ is a weak homotopy equivalence. 
\end{teo}

The process of obtaining $X/U_x$ from $X$ is called an \textit{open reduction} and the process of obtaining $X/\overline{\{x\}}$ from $X$ is called a \textit{closed reduction}. 

Osaki asserts in \cite{Osa} that he does not know whether by a sequence of reductions, each finite $T_0$-space can be reduced to 
the smallest space with the same homotopy groups.

\medskip

We show with the following  example that the answer to this question is negative.

\medskip

Let $X=\{a_1,b,a_2,c,d,e\}$ be the $6$-point $T_0$-space with the following order: $c,d<a_1; \ c,d,e<b$ and $d,e<a_2$. Let $D_3=\{c,d,e\}$ be the $3$-point discrete space and $Y=\mathbb{S}D_3=\{a,b,c,d,e\}$ the non-Haussdorf suspension of $D_3$. 

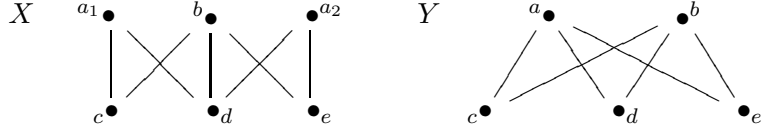
\begin{displaymath}
\xymatrix@C=10pt{ X & ^{a_1} \bullet \ \ \ \ar@{-}[d] \ar@{-}[rd] & ^{b} \bullet \ \ar@{-}[dl] \ar@{-}[d] \ar@{-}[rd] & \ \ \ \bullet ^{a_2} \ar@{-}[dl] \ar@{-}[d] \\ 
		& _c \bullet \ & \ \ \bullet _d & \ \ \bullet _e } 
\qquad
\xymatrix@C=4pt{ Y & & ^a \bullet \ \ar@{-}[dl] \ar@{-}[dr] \ar@{-}[drrr] & & \ \bullet ^b \ar@{-}[dlll] \ar@{-}[dl] \ar@{-}[dr]\\ 
		& _c \bullet \ & & \ \ \bullet _d & & \ \bullet _e } 
\end{displaymath}

The function $f:X\rightarrow Y$ defined by $f(a_1)=f(a_2)=a$, $f(b)=b$, $f(c)=c$, $f(d)=d$ and $f(e)=e$ is continuous because it preserves the order.

In order to prove that $f$ is a weak homotopy equivalence we  use Theorem $6$ of \cite{Mcc}.

\medskip

The sets $U_y$ form a basis-like covering of $Y$. 
Using the theory developed by Stong, it is  easy to verify that $f^{-1}(U_y)$ is contractible for each $y\in Y$ 
and, since $U_y$ is also contractible, the map $f|_{f^{-1}(U_y)}:f^{-1}(U_y)\rightarrow U_y$ is a weak homotopy equivalence for each $y\in Y$. 

Applying Theorem $6$ of \cite{Mcc}, one proves that $f$ is a weak homotopy equivalence. Therefore $X$ and $Y$ have the same homotopy groups. 

Another way to show that $X$ and $Y$ are weakly equivalent is considering the associated polyhedra $|\k (X)|$ and $|\k (Y)|$ which are 
homotopy equivalent to $S^1\vee S^1$.\\

On the other hand, it is easy to see that Osaki reduction methods cannot be applied to the space $X$. Therefore 
his methods are not effective in this case since we cannot obtain, by a sequence of reductions, the smallest space with the same homotopy groups of $X$.\\

\medskip

In order to achieve our goal, we must then choose a different approach.\\

We denote by $h(X)$ the height of a poset $X$, i.e. the maximum length of a chain in $X$.


\begin{teo} \label{esfera1}
Let $X\neq *$ be a minimal finite space. Then $X$ has at least $2h(X)$ points. 
Moreover, if $X$ has exactly $2h(X)$ points, then it is homeomorphic to $\mathbb{S}^{h(X)-1}S^0$.
\end{teo}

\begin{proof}
Let $x_1<x_2<\ldots <x_h$ be a chain in $X$ of length $h=h(X)$. Since $X$ is a minimal finite space, $x_i$ is not an up beat point for any $1\le i<h$. 
Then, for every $1\le i<h$ there exists $y_{i+1}\in X$ such that $y_{i+1}>x_i$ and $y_{i+1}\ngeq x_{i+1}$. 
We assert that the points $y_i$ (for $1<i\le h$) are all distinct from each other and also different from the $x_j$ ( $1\le j\le h$). 

\medskip

Since $y_{i+1}>x_i$, it follows that $y_{i+1}\neq x_j$ for all $j\le i$. But $y_{i+1}\neq x_j$ for all $j>i$ because $y_{i+1}\ngeq x_{i+1}$.

If $y_{i+1}=y_{j+1}$ for some $i<j$, then $y_{i+1}=y_{j+1}\ge x_j\ge x_{i+1}$, which is a contradiction.

\medskip

Since finite spaces with minimum or maximum are contractible and  $X\neq *$ is a minimal finite space, it cannot have a minimum. 
Then there exists $y_1\in X$ such that $y_1\ngeq x_1$. Therefore, $y_1$ must be distinct from the other $2h-1$ points and $\# X\ge 2h$.\\

Let us suppose now that $X$ has exactly $2h$ points, i.e.  $$X=\{x_1,x_2,\ldots ,x_h,y_1,y_2,\ldots ,y_h \}.$$

Because of the maximality of the chain $x_1<\ldots <x_h$, we get that $x_i$ and $y_i$ are incomparable for all $i$.

We show that $y_i<x_j$ and $y_i<y_j$ for all $i<j$ by induction in $j$.

\smallskip

For $j=1$ there is nothing to prove.
    
Let $1\le k <h$ and assume the statement holds for $j=k$. As $x_{k+1}$ is not a down beat point, there exists $z\in X$ such that $z<x_{k+1}$, and $z\nleq x_k$. 
Since $x_{k+1}$ and $y_{k+1}$ are incomparable, it follows that $z\neq y_{k+1}$. By induction we know 
that every point in $X$, with the exception of $y_k$ and $y_{k+1}$, is greater than $x_{k+1}$ or less than $x_k$. Then $z=y_k$ and so, $y_k<x_{k+1}$. 

Analogously, $y_{k+1}$ is not a down beat point and there exists $w\in X$ such that $w<y_{k+1}$ and $w\nleq x_k$. 
Again by induction, and because $y_{k+1}\ngeq x_{k+1}$, we deduce that $w$ must be $y_k$ and then $y_k<y_{k+1}$.

Furthermore, if $i<k$, then $y_i<x_k<x_{k+1}$ and $y_i<x_k<y_{k+1}$.\\

We proved that, for any $i<j$, we have that 
 $y_i<x_j$, $y_i<y_j$, $x_i<x_j$ and $x_i<y_j$. Moreover, for any $1\le i\le h$, 
$x_i$ and $y_i$ are incomparable.

This is exactly the order of $\mathbb{S}^{h-1}S^0$. Therefore $X$ is homeomorphic to $\mathbb{S}^{h-1}S^0$.
\end{proof}

\begin{teo} \label{esfera2} Any space with  the same homotopy groups of $S^n$ has at least $2n+2$ points. Moreover, $\S^n S^0$ is the unique space with $2n+2$ points with
this
property.
\end{teo}

\begin{proof}
The case $n= 1$ is trivial.

In other case, let us suppose that $X$ is a finite space with minimum cardinality such that $\pi_k(X,x)=\pi_k(S^n,s)$ for all $k\ge 0$. 
Then $X$ must be a minimal finite space and so $T_0$.

 By Hurewicz Theorem, $H_n(|\k(X)|)=\pi_n(|\k(X)|)=\pi_n(S^n) \neq 0$. This implies that  
the dimension of the simplicial complex $\k(X)$ must be at least $n$, which  means that the height of $X$ is at least $n+1$.

The result now follows immediately from the previous theorem.
\end{proof}

\begin{coro}
The $n$-sphere has a unique minimal finite model and it has $2n+2$ points.
\end{coro}

\begin{obs}
After concluding this paper, we found an old article of McCord (\textit{Singular homology and homotopy groups of finite spaces}, Notices of the American Mathematical Society, vol. 12(1965)) with a result (Theorem 2) without proof, from which the first part of \ref{esfera2} could be deduced. 
McCord's result can be easily deduced from our stronger theorem \ref{esfera1} (which also implies the uniqueness of these minimal models). 

Furthermore, we think that the proof of \ref{esfera1} itself is interesting because it relates the combinatorial methods of Stong's theory 
with McCord's point of view.  
\end{obs}

%
%
%
%

\section{Loops in the Hasse diagram and the fundamental group}

In this section we give a full description of the fundamental group of a finite $T_0$-space 
in terms of its Hasse diagram. This characterization is induced from the well known description of the fundamental group of a simplicial 
complex \cite{Spa}.

\smallskip

\begin{defi}
Let $(X,x_0)$ be a finite pointed $T_0$-space. An ordered pair of points $e=(x,y)$ is called an $\mathcal{H}$\textit{-edge of $X$} 
if $(x,y)\in \mathtt{E}(\mathcal{H}(X))$ or $(y,x)\in \mathtt{E}(\mathcal{H}(X))$. The point $x$ is called the \textit{origin} of $e$ and denoted $x=\o(e)$, the point $y$ is 
called the \textit{end} of $e$ and denoted $y=\e(e)$. The \textit{inverse} of an $\mathcal{H}$-edge $e=(x,y)$ is the $\mathcal{H}$-edge $e^{-1}=(y,x)$. 

An $\mathcal{H}$\textit{-path} in $(X,x_0)$ is a finite sequence of $\mathcal{H}$-edges (maybe empty) $\xi =e_1e_2\ldots e_n$ such that $\e(e_i)=\o(e_{i+1})$ for all $1\le i\le n-1$. The \textit{origin} of a non empty $\mathcal{H}$-path $\xi$ is $\o(\xi)=\o(e_1)$ and its \textit{end} is $\e(\xi)=\e(e_n)$. The origin and the end of the empty $\mathcal{H}$-path is $\o (\emptyset)=\e (\emptyset)=x_0$. If $\xi=e_1e_2 \ldots e_n$, we define $\overline{\xi}=e_n^{-1}e_{n-1}^{-1} \ldots e_1^{-1}$. If $\xi, \xi'$ are $\mathcal{H}$-paths such that $\e(\xi)=\o(\xi')$, we define the product $\mathcal{H}$-path $\xi \xi'$ as the concatenation of the sequence $\xi$ followed by the sequence $\xi'$.

An $\mathcal{H}$-path $\xi=e_1e_2\ldots e_n$ is said to be \textit{monotonic} if $e_i\in \mathtt{E}(\mathcal{H}(X))$ for all $1\le i\le n$ or $e_i^{-1}\in \mathtt{E}(\mathcal{H}(X))$ for all $1\le i\le n$.

A \textit{loop at $x_0$} is an $\mathcal{H}$-path that starts and ends in $x_0$. Given two loops $\xi,\xi'$ at $x_0$, we say that 
they are \textit{close} if there exist $\mathcal{H}$-paths $\xi_1,\xi_2,\xi_3,\xi_4$ such that $\xi_2$ and $\xi_3$ are monotonic and the set $\{\xi,\xi'\}$ coincides with $\{\xi_1\xi_2\xi_3\xi_4,\xi_1\xi_4\}$.

We say that two loops $\xi,\xi'$ at $x_0$  are $\mathcal{H}$\textit{-equivalent} if there exist a finite sequence 
of loops $\xi=\xi_1,\xi_2,\ldots ,\xi_n=\xi'$ such that any two consecutive are close. We denote by $\langle \xi \rangle$ the $\mathcal{H}$-equivalence 
class of a loop $\xi$ and $\mathscr{H}(X,x_0)$ the set of these classes.

\end{defi}

\begin{teo}
Let $(X,x_0)$ be a pointed  finite $T_0$-space. Then the product $\langle \xi \rangle \langle \xi' \rangle =\langle \xi \xi' \rangle$ 
is well defined and induces a group structure on $\mathscr{H}(X,x_0)$.
\end{teo}
\begin{proof}
It is easy to check that the product is well defined, associative and that  $\langle \emptyset \rangle$ is the identity. 
In order to prove that the inverse of $\langle e_1e_2\ldots e_n \rangle$ is $\langle e_n^{-1}e_{n-1}^{-1}\ldots e_1^{-1}\rangle$  we need 
to show that for any composable $\mathcal{H}$-paths $\xi, \xi'$ such that $\o (\xi)=\e (\xi')=x_0$ and for any $\mathcal{H}$-edge $e$, composable with $\xi$, one has that $\langle \xi e e^{-1}\xi'\rangle=\langle \xi \xi'\rangle$. 
But this follows immediately from the definition of close loops since $e$ and $e^{-1}$ are monotonic. 
\end{proof}

\begin{teo}
Let $(X,x_0)$ be a pointed finite $T_0$-space. Then  the edge-path group $E(\k(X),x_0)$ of $\k (X)$ with base vertex $x_0$ 
is isomorphic to $\mathscr{H}(X,x_0)$.
\end{teo}
\begin{proof}
Let us define 
\begin{align*}
\varphi :\mathscr{H}(X,x_0) & \longrightarrow E(\k(X),x_0) \\
\langle e_1e_2\ldots e_n\rangle & \longmapsto [e_1e_2\ldots e_n] \\
\langle \emptyset \rangle & \longmapsto [(x_0,x_0)]
\end{align*}

Where $[\xi]$ denotes the class of $\xi$ in $E(\k(X),x_0)$.\\

To prove that $\varphi$ is  well defined, let us suppose 
that the loops $\xi_1 \xi_2 \xi_3 \xi_4$ and $\xi_1 \xi_4$ are close, where $\xi_2=e_1e_2\ldots e_n$, $\xi_3=e_1'e_2'\ldots e_m'$ are monotonic $\mathcal{H}$-paths. 
By induction,  it can be proved 
that $[\xi_1 \xi_2 \xi_3 \xi_4]=[\xi_1 e_1e_2\ldots e_{n-j} (\o(e_{n-j+1}),\e(e_n) ) \xi_3 \xi_4]$ for $1\le j\le n$. 
In particular $[\xi_1 \xi_2 \xi_3 \xi_4]=[\xi_1 (\e(\xi_1), \e(e_n)) \xi_3 \xi_4]$. 

\medskip

Analogously, $$[\xi_1 (\e(\xi_1), \e(e_n)) \xi_3 \xi_4]=[\xi_1 (\e(\xi_1), \e(e_n)) (\o(e_1'), \o(\xi_4)) \xi_4]$$
 and then $$[\xi_1 \xi_2 \xi_3 \xi_4]=[\xi_1 (\e(\xi_1), \e(e_n)) (\o(e_1'), \o(\xi_4)) \xi_4]=[\xi_1 (\e(\xi_1), \e(e_n)) (\e(e_n), \e(\xi_1)) \xi_4]=$$ $$=[\xi_1 (\e(\xi_1), \e(\xi_1)) \xi_4]=[\xi_1 \xi_4]$$

\medskip

If $\xi=(x_0,x_1)(x_1,x_2)\ldots (x_{n-1},x_n)$ is an edge path in $\k (X)$ with $x_n=x_0$, then $x_{i-1}$ and $x_i$ are comparable for all $1\le i \le n$. In this case, we can find monotonic $\mathcal{H}$-paths $\xi_1, \xi_2, \ldots ,\xi_n$ such that $\o(\xi_i)=x_{i-1}, \ \e(\xi_i)=x_i$ for all $1\le i \le n$. Let us define 
\begin{align*}
\psi: E(\k(X),x_0) & \longrightarrow \mathscr{H}(X,x_0) \\
[\xi] & \longmapsto \langle \xi_1 \xi_2 \ldots \xi_n\rangle
\end{align*}

This definition does not depend on the election of the $\mathcal{H}$-paths $\xi_i$ since if two elections differ only for $i=k$ then $\xi_1 \ldots \xi_k \ldots \xi_n$ and $\xi_1 \ldots \xi_k' \ldots \xi_n$ are $\mathcal{H}$-equivalent because both of them are close to $\xi_1 \ldots \xi_k \xi_k^{-1} \xi_k' \ldots \xi_n$.

\medskip

The definition of $\psi$ does not depend on the representative. 
Let us suppose that $\xi' (x,y) (y,z) \xi''$ and $\xi' (x,z) \xi''$ are simply equivalent edge paths in $\k(X)$ that start and end in $x_0$, where $\xi$ and $\xi'$ are edge paths and $x,y,z$ are comparable. 

In the case that $y$ lies between $x$ and $z$, we can choose the monotonic $\mathcal{H}$-path corresponding to $(x,z)$ to be the yuxtaposition of the corresponding to $(x,y)$ and $(y,z)$. And so $\psi$ is equally defined in both edge paths.

In the case that $z\le x\le y$ we can choose monotonic $\mathcal{H}$-paths $\alpha$, $\beta$ from $x$ to $y$ and from $z$ to $x$. And then $\alpha$ will be the corresponding $\mathcal{H}$-path to $(x,y)$, $\overline{\alpha} \overline{\beta}$ the corresponding to $(y,z)$ and $\overline{\beta}$ to $(x,z)$. It only remains to prove that $\langle \gamma' \alpha \overline{\alpha} \overline{\beta} \gamma'' \rangle=\langle \gamma' \overline{\beta} \gamma'' \rangle$ for $\mathcal{H}$-paths $\gamma'$ and $\gamma''$, which is trivial.

The other cases are analogous to the last one.\\ 

It remains to verify that $\varphi$ and $\psi$ are mutually inverses, but this is clear. 
\end{proof}

Since $E(\k(X),x_0)$ is isomorphic to $\pi_1(|\k(X)|,x_0)$ (cf. \cite{Spa}),  we obtain the following result. 

\begin{coro}
Let $(X,x_0)$ be a pointed finite $T_0$-space. Then $\mathscr{H}(X,x_0)=\pi_1(X,x_0)$.
\end{coro}

\begin{obsi}
Since every finite space is homotopy equivalent to a finite $T_0$-space, this computation of the fundamental group can be applied to 
any finite space.
\end{obsi}

\bigskip

We finish this section with a couple of remarks on the Euler characteristic of finite spaces.

\begin{obs}
If $X$ is a finite $T_0$-space and $x<y\in X$ are such that $z>x$ implies $z\ge y$ ($x$ is an up beat point), then for any $z\in X$ 
such that $z$ is comparable with $x$, one has that  $z$ is also comparable with $y$.

An analogous proposition holds if $x>y$ are such that $z<x$ implies $z\le y$ ($x$ is a down beat point).
\end{obs}

Moreover, it is not difficult to prove the  following characterization of minimal finite spaces.

\begin{prop}
Let $X$ be a finite $T_0$-space. Then $X$ is a minimal finite space if and only if there are no $x,y\in X$ with $x\neq y$ such that if $z\in X$ is comparable with $x$, then so is with $y$.
\end{prop}

Since any finite $T_0$-space $X$ is weakly equivalent to the realization of $\k (X)$, whose simplices are the non empty chains in $X$, the Euler characteristic of $X$ is $$\chi (X)=\sum\limits_{C\in \mathcal{C}(X)} (-1)^{\# C +1} $$ where $\mathcal{C}(X)$ is the set of non empty chains of $X$. 

\medskip

Although it is very well known that the Euler characteristic is a homotopy invariant, 
we exhibit a basic proof of this fact in the case of finite spaces: 

\begin{teo}
Let $X$ and $Y$ be finite $T_0$-spaces with the same homotopy type. Then $\chi (X)=\chi (Y)$.
\end{teo}
\begin{proof}
Following \cite{Sto}, there exist two secuences of finite $T_0$-spaces 
$X=X_0\supseteq \ldots \supseteq X_n=X_c$ and $Y=Y_0\supseteq \ldots \supseteq Y_m=Y_c$,  where $X_{i+1}$ is constructed 
 from $X_i$ by removing a beat point and $Y_{i+1}$ is constructed from $Y_i$. 

Since $X$ and $Y$ are homotopy equivalent, $X_c$ and $Y_c$ are homeomorphic. Thus, $\chi (X_c)=\chi (Y_c)$. 

It suffices to show that the Euler characteristic does not change when a beat point is removed.

Let $P$ be a finite poset and let $p\in P$ be a beat point. 
By the result of above, there exists $q\in P$ such that $r$ comparable with $p$ implies $r$ comparable with $q$.

Hence we have a bijection

\begin{align*}
\noindent \varphi: \{C\in \mathcal{C}P \ | \ p\in C, \ q\notin C\} & \longrightarrow \{C\in \mathcal{C}P \ | \ p\in C, \ q\in C\} \\
C & \longmapsto C\cup \{q\}
\end{align*}

Therefore $$\chi (P)-\chi (P\smallsetminus\{p\})=\underset{p\in C\in \mathcal{C}P}{\sum}(-1)^{\# C+1}=\underset{q\notin C\ni p}{\sum}(-1)^{\# C+1}+\underset{q\in C\ni p}{\sum}(-1)^{\# C+1}=$$ $$=\underset{q\notin C\ni p}{\sum}(-1)^{\# C+1}+\underset{q\notin C\ni p}{\sum}(-1)^{\# \varphi (C) +1}=\underset{q\notin C\ni p}{\sum}(-1)^{\# C +1}+\underset{q\notin C\ni p}{\sum}(-1)^{\# C}=0.$$
\end{proof}


\section{Minimal finite models of graphs}

\begin{obs} \label{eulerdicetodo}
If $X$ is a connected finite $T_0$-space of height two, then $|\k(X)|$ is a connected graph, i.e. a CW complex of dimension one. Therefore, 
 the weak homotopy type of $X$ is completely determined by its Euler characteristic. 
More precisely, if $\chi (X)=\#X - \# \mathtt{E}(\mathcal{H}(X))=n$, then $X$ is a finite model of $\bigvee\limits_{i=1}^{1-n}S^1$.
\end{obs}

\begin{prop}
Let $X$ be a connected finite $T_0$-space and let $x_0\neq x\in X$ such that $x$ is neither maximal nor minimal in $X$. 
Then the inclusion map of the associated simplicial complexes $\k (X\smallsetminus \{x\})\subseteq \k (X)$ induces an epimorphism 
 $i_*:E(\k(X\smallsetminus \{x\}),x_0)\to E(\k(X),x_0)$ between their edge-path groups. 
\end{prop}

\begin{proof}
We  have to check that every closed edge path in $\k (X)$ with base point $x_0$ is equivalent to another edge path that does not go through $x$.

Let us suppose that $y\le x$ and $(y,x)(x,z)$ is an edge path in $\k (X)$.

 If $x\le z$ then $(y,x)(x,z)\equiv (y,z)$. 
In the case that $z<x$, since $x$ is not maximal in $X$, there exists $w>x$. Therefore
 $(y,x)(x,z)\equiv (y,x)(x,w)(w,x)(x,z)\equiv (y,w)(w,z)$. 

The case  $y\ge x$ is analogous.

This way one can eliminate $x$ from the writing of any closed edge path with base point $x_0$.
\end{proof}

It is important to note that the space $X\smallsetminus \{x\}$ of the previous proposition is also connected.

\medskip

As we pointed out in the introduction, the above result shows how powerful the combinatorial tools from finite spaces can be.
The conditions of maximality or minimality of points in a finite space  are
hard to express in terms of simplicial complexes. 

\begin{obs}
If $X$ is a finite $T_0$-space, then $h(X)\le 2$ if and only if every point in $X$ is maximal or minimal.
\end{obs}

\begin{coro}
Let $X$ be a connected finite space. 
Then there exists a connected  $T_0$-subspace $Y\subseteq X$ of height at most two such that the fundamental group of $X$ is a 
quotient of the fundamental group of $Y$.
\end{coro}
\begin{proof}
We can assume that $X$ is $T_0$ because $X$ has a core. Since the edge-path group is isomorphic to the fundamental group, the result follows
straight forward from the previous proposition. 
\end{proof}

\begin{obs}
Note that the fundamental group of a connected finite $T_0$-space of height at most two is finitely generated by \ref{eulerdicetodo}. Therefore, the path-connected spaces whose fundamental group does not have a finite set of generators do not admit finite models.
\end{obs}

\begin{coro}
Let $n\in \mathbb{N}$. If $X$ is a minimal finite model of $\bigvee\limits_{i=1}^{n}S^1$, then $h(X)=2$.
\end{coro}

\begin{proof}
Let $X$ be a minimal finite model of $\bigvee\limits_{i=1}^{n}S^1$. Then there exists a connected $T_0$-subspace $Y\subseteq X$ of height two, $x\in Y$ and an epimorphism 
from $\pi_1(Y,x)$ to $\pi_1(X,x)=\overset{n}{\underset{i=1}{\ast}}\mathbb{Z}$.

Since $h (Y)=2$, $Y$ is a model of a graph, thus $\pi_1(Y,x)=\overset{m}{\underset{i=1}{\ast}}\mathbb{Z}$ for some integer $m$.

 Note that $m \geq n$.
 
 \medskip
 
There are $m$ edges of $\mathcal{H}(Y)$ which are not in a maximal tree of the underlying 
non directed graph of $\mathcal{H}(Y)$ (i.e. $\k (Y)$). Therefore, we can remove $m-n$ edges from $\mathcal{H}(Y)$ in such a way that it remains connected and the new space $Z$ obtained this way is a model of $\bigvee\limits_{i=1}^{n}S^1$.

Note that $\# Z=\# Y\le \# X$. Since $X$ is a minimal finite model, then $\#X\le \#Z$. Therefore $X=Y$ has height two.
\end{proof}

If $X$ is a minimal finite model of $\bigvee\limits_{i=1}^{n}S^1$ and we call $i=\# \{y\in X \ | \ y$ is maximal$\}$, $j=\# \{y\in X \ | \ y$ is minimal$\}$, then $\# X=i+j$ and $\# \mathtt{E}(\mathcal{H}(X))\le ij$. Since $\chi (X)=1-n$, we have that $n\le ij-(i+j)+1=(i-1)(j-1)$.\\

We can now state the main result of this section.

\begin{teo}
Let $n\in \mathbb{N}$. A finite $T_0$-space $X$ is a minimal finite model of $\bigvee\limits_{i=1}^{n}S^1$ if and only if $h(X)=2$, $\# X=min\{ i+j \ | \ (i-1)(j-1)\ge n\}$ and $\# \mathtt{E}(\mathcal{H}(X))=\# X +n-1$.
\end{teo}
\begin{proof}
We have already proved that if $X$ is a minimal finite model of $\bigvee\limits_{i=1}^{n}S^1$, then $h(X)=2$ and $\#X \ge min\{ i+j \ | \ (i-1)(j-1)\ge n\}$. 

If $i$ and $j$ are such that $n\le (i-1)(j-1)$, we can consider $Y=\{x_1,x_2,\ldots ,x_i,y_1,y_2,\ldots y_j\}$ with the order $y_k\le x_l$ for all $k,l$, which is a model of $\bigvee\limits_{k=1}^{(i-1)(j-1)}S^1$. Then we can remove $(i-1)(j-1)-n$ edges from $\mathcal{H} (X)$ to obtain a connected space of cardinality $i+j$ which is a finite model of $\bigvee\limits_{k=1}^{n}S^1$. Therefore $\#X\le \#Y=i+j$. 

This is true for any $i,j$ with $n\le (i-1)(j-1)$, then $\#X=min\{ i+j \ | \ (i-1)(j-1)\ge n\}$.

Moreover, $\# \mathtt{E}(\mathcal{H}(X))=\# X +n-1$ because $\chi (X)=1-n$.\\

In order to show the converse of the theorem we only need to prove that the conditions 
$h(X)=2$, $\# X=min\{ i+j \ | \ (i-1)(j-1)\ge n\}$ and $\# \mathtt{E}(\mathcal{H}(X))=\# X +n-1$ imply that $X$ is connected, because in this case, by \ref{eulerdicetodo}, the first and third conditions would say that $X$ is a model of $\bigvee\limits_{i=1}^{n}S^1$, and the second condition would say that it has the right cardinality.

\medskip

Suppose $X$ satisfies the conditions of above and let  $X_l$,  $1\le l \le k$, be the connected components of $X$. 

Let us denote by $M_l$ the set of maximal elements of $X_l$ and let $m_l=X_l\smallsetminus M_l$. 
Let $i=\sum\limits_{r=1}^k \#M_l$, $j=\sum\limits_{r=1}^k \#m_l$. 

Since $i+j=\# X=min\{ s+t \ | \ (s-1)(t-1)\ge n\}$, it follows that $(i-2)(j-1)<n=\# \mathtt{E}(\mathcal{H}(X))-\# X +1=\# \mathtt{E}(\mathcal{H}(X))-(i+j) +1$. 
Hence $ij-\# \mathtt{E}(\mathcal{H}(X))<j-1$.\\ 

This means that $\k (X)$ differs from the complete bipartite graph $(\cup m_l,\cup M_l)$ in less than $j-1$ edges. 

Since there are no edges from $m_r$ to $M_l$ if $r\neq l$, $$j-1>\sum\limits_{l=1}^k \# M_l(j-\#m_l)\ge \sum\limits_{l=1}^k (j-\#m_l)=(k-1)j$$ 

Therefore $k=1$ and the proof is complete.
\end{proof}

\begin{obs}
The cardinality of a minimal finite model of $\bigvee\limits_{i=1}^{n}S^1$ is 
$$min \{2\lceil \sqrt{n}+1 \rceil, 2\bigg\lceil  \frac{1+\sqrt{1+4n}}{2} \bigg\rceil +1 \}.$$
\end{obs}

Note that a space may admit many minimal finite models as we can see in the following example.

\begin{ej}
Any  minimal finite model of $\bigvee\limits_{i=1}^{3}S^1$ has $6$ points and $8$ edges. So, they are, up to homeomorphism

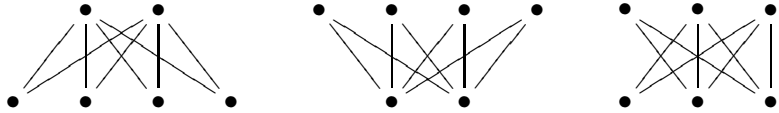
\begin{displaymath}
\xymatrix@C=16pt{ & \bullet \ar@{-}[ld] \ar@{-}[d] \ar@{-}[rd] \ar@{-}[rrd] & \bullet \ar@{-}[ld] \ar@{-}[d] \ar@{-}[rd] \ar@{-}[lld]  & \\
		\bullet & \bullet & \bullet & \bullet }
\qquad
\xymatrix@C=16pt{ \bullet & \bullet & \bullet & \bullet \\
		& \bullet \ar@{-}[lu] \ar@{-}[u] \ar@{-}[ru] \ar@{-}[rru] & \bullet \ar@{-}[lu] \ar@{-}[u] \ar@{-}[ru] \ar@{-}[llu]  & }
\qquad
\xymatrix@C=16pt{ \bullet \ar@{-}[rd] \ar@{-}[rrd] & \bullet \ar@{-}[ld] \ar@{-}[d] \ar@{-}[rd] & \bullet \ar@{-}[lld] \ar@{-}[ld] \ar@{-}[d] \\
		\bullet & \bullet & \bullet }
\end{displaymath}

\end{ej}

\medskip

In fact, it is not hard to prove, using our characterization, that $\bigvee\limits_{i=1}^{n}S^1$ has a unique minimal finite model if and only if $n$ is a square.

\medskip

\begin{obsi}
Since any graph is a $K(G,1)$, the minimal finite models of a graph $X$ are, in fact, the smallest spaces with the same homotopy groups of $X$.
\end{obsi}

\email{jbarmak@dm.uba.ar, gminian@dm.uba.ar}


\begin{thebibliography}{99}

\bibitem{Ale} P.S. Alexandroff. \textit{Diskrete R\"aume}.
    MathematiceskiiSbornik (N.S.) 2(1937), 501-518.
  


\bibitem{May} J.P. May. \textit{Finite topological spaces}.
    Notes for REU (2003).

\bibitem{May2} J.P. May. \textit{Finite spaces and simplicial complexes}.
    Notes for REU (2003).
    
\bibitem{May3} J.P. May. \textit{Finite groups and finite spaces}.
    Notes for REU (2003).

\bibitem{May4} J.P. May. \textit{ A concise course in algebraic topology}.  Chicago lecture notes in mathematics (1999).
   
\bibitem{Mcc} M.C. McCord. \textit{Singular homology groups and homotopy groups of finite topological spaces}.
    Duke Mathematical Journal 33(1966), 465-474.
    
\bibitem{Osa} T. Osaki. \textit{Reduction of finite topological spaces}.
    Interdiciplinary Information Sciences 5(1999), 149-155.  
  
\bibitem{Spa} E. Spanier. \textit{Algebraic
    Topology}. Springer (1966).
    
\bibitem{Sto} R.E. Stong. \textit{Finite topological spaces}.
    Trans. Amer. Math. Soc. 123(1966), 325-340.
    
\end{thebibliography}
\end{document}